\documentclass[a4paper,11pt]{amsart}

\usepackage {amssymb}
\usepackage{amscd}

\theoremstyle{plain}
\newtheorem{thm}{Theorem}[section]
\newtheorem{lemma}[thm]{Lemma}
\newtheorem{cor}[thm]{Corollary}
\newtheorem{prop}[thm]{Proposition}

\theoremstyle{definition}
\newtheorem{defi}[thm]{Definition}

\def\dim{\mathop{\hbox {dim}}\nolimits}

\newcommand{\frg}{\mathfrak{g}}

\newcommand{\frl}{\mathfrak{l}}

\newcommand{\frn}{\mathfrak{n}}

\newcommand{\frr}{\mathfrak{r}}
\newcommand{\frs}{\mathfrak{s}}

\newcommand{\frgl}{\mathfrak{gl}}
\newcommand{\frsl}{\mathfrak{sl}}

\newcommand{\bbC}{\mathbb{C}}

\begin{document}
\title{Refinement of Ado's Theorem in Low Dimensions and Application in Affine Geometry}

\author{Yi-Fang Kang}
\address{Department of Mathematics, Tianjin University,
Tianjin, 300072, PR China} \email{kangyf\char'100tju.edu.cn}

\author{Cheng-Ming Bai$^*$}
\address{Chern Institute of Mathematics \& LPMC, Nankai University,  Tianjin, 300071, PR China}
\email{baicm\char'100nankai.edu.cn}

\def\shorttitle{Refinement of Ado's Theorem and Affine Geometry} % 压缩标题, 尽量不超过60个字符 (包括空格)

\abstract {In this paper, we construct a faithful representation
with the lowest dimension for every complex Lie algebra in
dimension $\leq 4$. In particular, in our construction, in the
case that the faithful representation has the same dimension of
the Lie algebra, it can induce an \'etale affine representation
with base zero which has a natural and simple form and gives a
compatible left-symmetric algebra on the Lie algebra. Such affine
representations do not contain any nontrivial one-parameter
subgroups of translation.}
\endabstract

\thanks{*Corresponding author}

%\footnote{\textup{2000} Mathematics Subject Classification. 17B,
%53C}
\subjclass[2000]{17B, 53C}

\keywords{Lie algebra, Ado's Theorem, Left-symmetric algebra,
Affine structure}

\maketitle

\section{Introduction}

Let $\frg$ be an $n$-dimensional Lie algebra over a field ${\mathbb
F}$ of characteristic zero. A well known theorem due to Ado states
that there exists a finite-dimensional faithful representation of
$\frg$ ([J]). That is, there exists an injective homomorphism
$f:\frg\rightarrow \frg\frl(N)$ of Lie algebras for some $N\in
{\mathbb N}$. It is natural to ask what the value of $N$ is in terms
of the dimension of $\frg$.

Define an invariant of $\frg$ by ([Bu1])
$$\mu(\frg)={\rm min}\{ \dim_{\mathbb F}\rho|\rho\; {\rm is}\;\;{\rm
a}\;\; {\rm faithful}\;\;{\rm representation}\;\;{\rm
of}\;\;\frg\}.\eqno (1.1)$$ By Ado's Theorem, $\mu(\frg)$ is
finite. It is also known that $\mu (\frg)\leq f(n)$ for a function
$f$ only depending on $n$.

In general, except some special cases, it is hard to decide the
exact $\mu(\frg)$. Even it is not easy to give a suitable estimation
on $\mu (\frg)$, either, especially in the cases of solvable and
nilpotent Lie algebras. Up to now, most of the estimation on $\mu$
for a general nilpotent or solvable Lie algebra is quite ``big". For
example, for an $n$-dimensional complex solvable Lie algebra $\frg$,
there is an estimation $\mu(\frg)<1+n+n^n$ ([R]). And when $\frg$ is
an $n$-dimensional nilpotent Lie algebra with nilpotency class $k$
(that is, $\frg^{k}\ne 0, \frg^{k+1}=0$, where
$\frg^{i}=[\frg,\frg^{i-1}]$, $\frg^1=\frg$), there is the bound
$$\mu(\frg)\leq \sum_{j=0}^k \binom{n-j}{k-j} p(j)<\frac{3}{\sqrt{n}}2^n,\eqno (1.2)$$
where $\binom{n-j}{k-j}$ is the binomial number and $p(j)$ is the
number of partitions of $j$ ([Bu1-2]).

On the other hand,  it is obvious that any $n$-dimensional Lie
algebra $\frg$ with trivial center satisfies $\mu(\frg)\leq n$. And
any $n$-dimensional real Lie algebra $\frg$ which has a codimension
one abelian ideal also satisfies $\mu(\frg)\leq n$ ([GST]).
Furthermore, a useful estimation on $\mu$ is given for a Lie algebra
$\frg$ with a compatible left-symmetric algebra structure which
corresponds to the affine structure on a Lie group $G$ whose Lie
algebra is $\frg$ ([M], [V]). For such a $n$-dimensional Lie algebra
$\frg$, it is known that $\mu(\frg)\leq n+1$. Left-symmetric
algebras are a kind of natural (nonassociative) algebraic systems
appearing in many fields in mathematics and mathematical physics,
such as affine and symplectic geometry ([K],[DM1-2]), integrable
systems ([Bo]), classical and quantum Yang-Baxter equation
([ES],[GS]), Poisson brackets and infinite-dimensional Lie algebras
([BN],[Z]), quantum field theory ([CK]), operads ([CL]) and so on.
Fortunately, there exist left-symmetric algebras on a lot of Lie
algebras, in particular, for most of the solvable and nilpotent Lie
algebras in low dimensions.

Therefore, it is practicable (and necessary) to decide the exact
$\mu$ for certain concrete Lie algebras, at least in low dimensions.
In this paper, we find all $\mu(\frg)$ for the complex Lie algebras
in dimension $\leq 4$. We also give the corresponding faithful
representations. The construction is nontrivial since in the case
that the faithful representation has the same dimension of the Lie
algebra, it can induce an \'etale affine representation with base
zero which has a natural and simple form and gives a compatible
left-symmetric algebra structure on the Lie algebra. Although it is
based on the example-computation, it can be regarded as a guide for
further development. We would like to point out that these results
are partly consistent with some results in [GST] ($\mu(\frg)$ for
the real indecomposable algebras are computed there). The paper is
organized as follows. In section 2, we give some preliminaries and
basic results. We omit most of the proofs since they can be found in
some related references.  In section 3, we get $\mu(\frg)$ for all
complex Lie algebras in dimension $\leq 4$ and construct the
corresponding faithful representations. In section 4, we construct
the induced \'etale affine representations and then get the
corresponding left-symmetric algebras.

Throughout the paper, all algebras are finite dimensional and over
the complex field $\bbC$.
%and $\langle\ |\ \rangle$ stands for an
%algebra with a basis and nonzero products at each side of $|$.

\section{Preliminaries and basic results}

For some special cases,  $\mu(\frg)$ can be decided ([Bu1-2]).

\begin{lemma}
Let $\frg$ be an $n$-dimensional commutative Lie algebra. Then
$\mu(\frg)=[2\sqrt{n-1}]$, where $[x]$ is the least integer greater
or equal than $x$.
\end{lemma}

\begin{lemma}
Let $\frg$ be a 2-step nilpotent Lie algebra in dimension $n$ with
1-dimensional center. Then $n$ is odd and $\mu(\frg)=(n+3)/2$.
\end{lemma}

Recall that a filform Lie algebra is an $n$-dimensional nilpotent
Lie algebra with nilpotency class $n-1$.

\begin{lemma}
Let $\frg$ be a filform Lie algebra in dimension $n$. If the
commutator of  $\frg$ is abelian or $n<10$, then $\mu(\frg)=n$.
\end{lemma}

Obviously, there is the (unique) 1-dimensional Lie algebra $\bbC$.
In dimension 2,  there are exactly two (non-isomorphic) Lie
algebras: abelian Lie algebra $\bbC^2$ and $\frr_2(\bbC)$ with a
basis  $\{e_1,e_2\}$ satisfying $[e_1,e_2]=e_1$. The classification
of 3 and 4 dimensional complex Lie algebras was summarized in [BS]
(only the non-zero products are given).
\begin{lemma}
Every complex 3-dimensional Lie algebra is exactly isomorphic to
one Lie algebra of the list shown in Table 1.
\end{lemma}
\begin{table}[b]\caption{}
\begin{tabular}{|c|c|c|}
\hline $\frg$& $\mu (\frg)$ &Lie brackets\\
\hline $\bbC^3$ &3&  --\\
\hline $\frn_3(\bbC)$ & 3 & $[e_1,e_2]=e_3$ \\
\hline $\frr_2(\bbC)\oplus\bbC$ & 2& $[e_1,e_2]=e_1$\\
\hline $\frr_3(\bbC)$& 3 &$[e_1,e_2]=e_2,[e_1,e_3]=e_2+e_3$\\
\hline $\frr_{3,\lambda}(\bbC)$& 3&$[e_1,e_2]=e_2,[e_1,e_3]=\lambda e_3,\lambda\in\bbC^*,|\lambda|<1$ or $\lambda=e^{i\theta}, 0\leq \theta\leq \pi$\\
\hline $\frsl_2(\bbC)$&2& $[e_1,e_2]=e_3,[e_1,e_3]=-2e_2,[e_2,e_3]=2e_2$\\
\hline
\end{tabular}
\end{table}
\begin{lemma}
Every complex 4-dimensional Lie algebra is exactly isomorphic to
one Lie algebra of the list shown in Table 2.
\end{lemma}

{\small \begin{table}[t]\caption{}
\begin{tabular}{|c|c|c|}
\hline $\frg$&$\mu(\frg)$ &Lie brackets\\
\hline $\bbC^4$ &4 & --\\
\hline $\frn_3(\bbC)\oplus\bbC$ & 3&$[e_1,e_2]=e_3$ \\
\hline $\frr_2(\bbC)\oplus\bbC^2$&3& $[e_1,e_2]=e_1$\\
\hline $\frr_3(\bbC)\oplus\bbC$&3 &$[e_1,e_2]=e_2,[e_1,e_3]=e_2+e_3$\\
\hline $\frr_{3,\lambda}(\bbC)\oplus\bbC$&3& $[e_1,e_2]=e_2,[e_1,e_3]=\lambda e_3,\lambda\in\bbC^*,|\lambda|<1$ or $\lambda=e^{i\theta}, 0\leq \theta\leq \pi$\\
\hline $\frr_2(\bbC)\oplus\frr_2(\bbC)$& 3&$[e_1,e_2]=e_1,[e_3,e_4]=e_3$\\
\hline $\frsl_2(\bbC)\oplus\bbC$& 2&$[e_1,e_2]=e_3,[e_1,e_3]=-2e_1,[e_2,e_3]=2e_2$\\
\hline $\frn_4(\bbC)$& 4& $[e_1,e_2]=e_3,[e_1,e_3]=e_4$\\
\hline $\frg_1(\alpha)$& 4&$[e_1,e_2]=e_2,[e_1,e_3]=e_3,[e_1,e_4]=\alpha e_4,\alpha\in\bbC^*$\\
\hline $\frg_2(\alpha,\beta)$&4& $[e_1,e_2]=e_3,[e_1,e_3]=e_4,[e_1,e_4]=\alpha e_2-\beta e_3+e_4$,\\ & &$\alpha\in\bbC^*,\beta\in \bbC$ or $\alpha,\beta=0$\\
\hline $\frg_3(\alpha)$&4& $[e_1,e_2]=e_3,[e_1,e_3]=e_4,[e_1,e_4]=\alpha (e_2+e_3),\alpha\in\bbC^*$\\
\hline $\frg_4$& 4&$[e_1,e_2]=e_3,[e_1,e_3]=e_4,[e_1,e_4]=e_2$\\
\hline $\frg_5$& 4&$[e_1,e_2]=\frac{1}{3}e_2+e_3,[e_1,e_3]=\frac{1}{3}e_3,[e_1,e_4]=\frac{1}{3}e_4$\\
\hline $\frg_6$& 3&$[e_1,e_2]=e_2,[e_1,e_3]=e_3,[e_1,e_4]=2e_4,[e_2,e_3]=e_4$\\
\hline $\frg_7$& 3&$[e_1,e_2]=e_3,[e_1,e_3]=e_2,[e_2,e_3]=e_4$\\
\hline $\frg_8(\alpha)$& 3($\alpha\ne \frac{1}{4}$)&$[e_1,e_2]=e_3,[e_1,e_3]=-\alpha e_2+e_3,[e_1,e_4]=e_4,[e_2,e_3]=e_4,\alpha\in\bbC$\\
& 4($\alpha=\frac{1}{4}$) &\\
\hline
\end{tabular}
\end{table}}

On the other hand, we have ([Bai], [BM1-3], [Bu1-4])

\begin{defi}
Let $A$ be a vector space over a field ${\mathbb F}$ with a
bilinear product $(x,y)\rightarrow xy$. $A$ is called a
left-symmetric algebra if for any $x,y,z\in A$, the associator
$$(x,y,z)=(xy)z-x(yz)\eqno (2.1)$$
is symmetric in $x,y$, that is,
$$(x,y,z)=(y,x,z),\;\;{\rm or}\;\;{\rm
equivalently}\;\;(xy)z-x(yz)=(yx)z-y(xz).\eqno (2.2)$$
\end{defi}

\begin{prop} Let $A$ be a left-symmetric
algebra. For any $x,y\in A$, let $L:A\rightarrow gl(A)$ by
$L(x)=L_x$, where $L_x$ denotes the left multiplication operator,
that is, $L_x(y)=xy$.

(1) The commutator
$$[x,y]=xy-yx,\;\;\forall x,y\in A,\eqno (2.3)$$
defines a Lie algebra ${\frg}(A)$, which is called the sub-adjacent
Lie algebra of $A$ and $A$ is called the compatible left-symmetric
algebra structure on the Lie algebra ${\frg}(A)$.

(2) $L$ gives a regular representation of the Lie algebra $\frg
(A)$, that is,
$$[L_x,L_y]=L_{[x,y]},\;\;\forall x,y\in A. \eqno (2.4)$$
\end{prop}

Let ${\frg}$ be a Lie algebra. An affine representation of $\frg$
is a Lie algebra homomorphism $\Phi:\frg\rightarrow {\rm aff}(V)$,
where ${\rm aff}(V)=\{ \phi=\left( \begin{matrix} \varphi & v\cr 0
& 0\cr \end{matrix}\right)|\varphi\in gl(V),v\in V\}\subset
gl(V\oplus {\bf F})$ is the Lie algebra of the group of affine
transformations ${\rm Aff}(V)=\{ g=\left( \begin{matrix} A & v\cr
0 & 1\cr\end{matrix}\right)|A\in {\rm GL}(V),v\in V\}\subset {\rm
GL}(V\oplus{\bf F})$. In addition, $\Phi$ is called an \'etale
affine representation with base point $v$ if there exists $v\in V$
such that the mapping $ev_v:\frg\rightarrow V$ defined by
$ev_v(x)=\Phi(x)v$ for any $x\in\frg$ is an isomorphism.
 It is known that $\Phi=(\rho,q):\frg\rightarrow {\rm aff}(V)$ with
$\Phi(x)=\rho(x)+q(x)$ is an affine representation if and only if
$\rho:\frg\rightarrow gl(V)$ is a representation of $\frg$ and
$q:\frg\rightarrow V$ is a linear map satisfying
$$q[x,y]=\rho(x)q(y)-\rho(y)q(x),\forall x,y\in {\frg}.\eqno (2.5)$$

\begin{prop} Let ${\frg}$ be a Lie algebra.
Then there is a compatible left-symmetric algebra structure on
${\frg}$ if and only if $\frg$ has an \'etale affine
representation.
\end{prop}

In fact, let $\Phi=(\rho,q):\frg\rightarrow {\rm aff}(V)$ be an
\'etale affine representation of $\frg$, then
$$x*y=ev_v^{-1}[\rho(x) ev_v(y)],\;\;\forall x,y\in \frg, \eqno (2.6)$$
defines a left-symmetric algebra structure on $\frg$. Conversely,
for a left-symmetric algebra $A$, $\Phi=(L,id)$ is an \'etale affine
representation of $\frg(A)$ with base $0$, where $id$ is the
identity transformation on ${\frg}(A)$. Moreover, if $A$ is a
left-symmetric algebra, then $N(A)=\{ x\in A|L_x=0\}$ is an ideal of
$A$ which is called a kernel ideal.

\begin{prop} {\bf ([M])}
A left-symmetric algebra $A$ has a zero kernel ideal if and only
if its corresponding \'etale affine representation  does not
contain any nontrivial one-parameter subgroups of translation.
\end{prop}

\begin{thm}  Let ${\frg}$ be a Lie algebra. If there exists an
\'etale affine representation $\Phi=(\rho,q)$ of $\frg$ such that
$\rho$ is a faithful representation of $\frg$, then the
corresponding left-symmetric algebra has the zero kernel ideal.
Therefore $\Phi$ does not contain any nontrivial one-parameter
subgroups of translation.
\end{thm}
\noindent {\bf Proof}\quad Let $x\in N(\frg)$. By equation (2.6),
we have $ev_v^{-1}[\rho(x) ev_v(y)]=0$ for every $y\in \frg$.
Since $ev_v$ is a linear isomorphism, $\rho(x)z=0$ for every $z\in
\frg$. Then $x\in {\rm Ker}\rho$. Since $\rho$ is faithful,
$x=0$.\hfill $\Box$

\section{The lowest dimensional faithful representations of Lie algebras in dimension $\leq 4$}

Obviously, $\mu (\frg)=1$ if only if $\frg$ is isomorphic to
$\bbC$.  Furthermore, by Lemma 2.1, we have
\begin{cor} With the notations given in section 2,
$$\mu (\bbC)=1,\;\mu(\bbC^2)=2,\;\mu(\bbC^3)=3,\mu(\bbC^4)=4.$$
\end{cor}

By Lemma 2.2, we have

\begin{cor} With the notations given in section 2,
$\mu(\frn_3(\bbC))=3$.
\end{cor}

By Lemma 2.3, we have

\begin{cor} With the notations given in section 2,
$\mu(\frn_4(\bbC))=4$.
\end{cor}

\begin{lemma} Let $\frg$ be a Lie algebra in dimension $n\leq 4$. If $\frg$ is 2-solvable
(that is, $[\frg,\frg]$ is a commutative Lie algebra) and $\dim
[\frg,\frg]=n-1$, then $\mu(\frg)\geq n$.
\end{lemma}

{\bf Proof}\quad We assume that $\frg$ has a faithful representation
$(\rho, V)$ in dimension $m<n$. Since $\frg$ is solvable, by Lie's
theorem, there exists a basis $\{ e_1,\cdots,e_m\}$ of $V$ such that
$\rho(x)\in \frg\frl(m,\bbC)$ is upper-triangle for all $x\in \frg$.
Therefore, $\rho([\frg,\frg])$ is in the space of strictly
upper-triangle matrices whose dimension is $m(m-1)/2$. When $1<n\leq
3$, we have $m(m-1)/2<n-1$, which is a contradiction. When $n=4$, we
can suppose $m=3$ since $\dim[\frg,\frg]=3$. Then
$\rho([\frg,\frg)]$ is exactly the Lie algebra of the $3\times 3$
strictly upper-triangle matrices which is not commutative. It is
contradictive to the fact that $[\frg,\frg]$ is commutative.
Therefore, the conclusion holds.\hfill $\Box$

\begin{cor} With the notations given in section 2,

(n=2) $\mu(\frr_2(\bbC))=2$;

(n=3) $\mu (\frr_3(\bbC))=\mu(\frr_{3,\lambda}(\bbC))=3$;

(n=4) $\mu(\frg_1(\alpha))=\mu(\frg_2(\alpha,\beta) (\alpha\in
\bbC^*,\beta\in
\bbC))=\mu(\frg_3(\alpha))=\mu(\frg_4)=\mu(\frg_5)=4$.
\end{cor}

{\bf Proof}\quad Let $\frg$ be one of the above Lie algebras. Then
$\frg$ satisfies the conditions in Lemma 3.4 and its adjoint
representation is faithful since $\frg$ is centerless. Therefore
$\mu(\frg)=n$.\hfill $\Box$

%By Theorem 3.3, we know both $\mu (\frr_3(\bbC)$ and
%$\mu(\frr_{3,\lambda}(\bbC))$ is greater than 2. Since both of
%them have the trivial center, their adjoint representation is
%faithful. Therefore the conclusion follows.

\begin{lemma}
Let $\frg$ be a complex centerless Lie algebra. Then
$\mu(\frg)=\mu (\frg\oplus \bbC)$.
\end{lemma}

\noindent {\bf Proof}\quad Let $\mu(\frg)=n$ and $\rho:
\frg\rightarrow \frg\frl(n,\bbC)$ be a faithful representation.
Obviously, $\mu (\frg\oplus \bbC)\geq n$.  Since the identity
matrix $I_n\notin \rho (\frg)$, $\rho (\frg)\oplus \bbC \cdot I_n$
gives a faithful representation of $\frg\oplus \bbC$. Therefore
$\mu (\frg\oplus \bbC)=\mu(\frg)=n$. \hfill $\Box$

\begin{cor} With the notations given in section 2,
$$\mu (\frr_2(\bbC)\oplus \bbC)=2,\;\mu (\frr_3(\bbC)\oplus \bbC)=\mu(\frr_{3,\lambda}(\bbC)\oplus\bbC)=3.$$
\end{cor}

\begin{thm}
Let $\frg$ be a complex Lie algebra. If $\mu(\frg)=2$, then $\frg$
is isomorphic to one of the following algebras:
$$\bbC^2, \frr_2(\bbC), \frr_2(\bbC)\oplus \bbC, \frs\frl_2(\bbC),
\frs\frl_2(\bbC)\oplus\bbC.$$
\end{thm}

\noindent {\bf Proof}\quad If $\mu(\frg)=2$, then there exists an
injective homomorphism $f:\frg\rightarrow
\frg\frl(2,\bbC)=\frs\frl_2(\bbC)\oplus\bbC$ of Lie algebras. Hence
$\dim \frg\leq 4$ and $\dim \frg=4$ if and only if $\frg$ is
isomorphic to $\frg\frl(2,\bbC)$. $\mu (\frs\frl(2,\bbC))=2$ since
the fundamental representation of $\frs\frl(2,\bbC)$ is faithful.
Therefore we have found all $\mu (\frg)$ when $\frg$ is a complex
Lie algebra in dimension 2 or 3. Hence the conclusion holds.\hfill
$\Box$

There are still the following 4-dimensional Lie algebras whose
$\mu (\frg)$s have not been decided:
$$\frn_3(\bbC)\oplus\bbC,
\frr_2(\bbC)\oplus\bbC^2,\frr_2(\bbC)\oplus\frr_2(\bbC),\frg_2(0,0),\frg_6,
\frg_7, \frg_8(\alpha).$$ Let $\frg$ be one of the above Lie
algebras. Then $\frg$ is solvable and $\mu(\frg)\geq 3$. It is known
that there exists a left-symmetric algebra structure on $\frg$
([Bu1]). So $\mu (\frg)\leq 5$ (in fact, $\mu(\frg)\leq 4$ from
[GST] or the following discussion). If we can construct a
3-dimensional faithful representation, then $\mu (\frg)=3$.
Otherwise, we need to show that $\mu (\frg)> 3$ and if we can
construct a 4-dimensional faithful representation, then $\mu
(\frg)=4$.

\begin{prop}
With the notations given in section 2, if $\frg$ is isomorphic to
one of the following algebras: $\frn_3(\bbC)\oplus\bbC,
\frr_2(\bbC)\oplus\bbC^2,\frr_2(\bbC)\oplus\frr_2(\bbC), \frg_6,
\frg_7, \frg_8(\alpha) (\alpha\ne \frac{1}{4}$), then
$\mu(\frg)=3$.
\end{prop}

{\bf Proof}\quad We can construct a 3-dimensional faithful
representation for every above Lie algebra. The construction is
given in  Table 3. \hfill $\Box$

\begin{prop}
With the notations given in section 2, if $\frg$ is isomorphic to
$\frg_2(0,0)$ or $\frg_8(\frac{1}{4})$, then $\mu(\frg)=4$.
\end{prop}

{\bf Proof}\quad  Let $\frg$ be $\frg_2(0,0)$ or
$\frg_8(\frac{1}{4})$ and $\{e_1, e_2, e_3, e_4\}$ be a basis of
$\frg$ with the Lie brackets in Table 2. Suppose that $(\rho, V)$
is a 3-dimensional faithful representation of $\frg$. By Lie's
theorem, we can select a basis of $V$ such that
$\rho(x)\in\frgl(3,\bbC)$ is upper-triangle for all $x\in\frg$.
That is, we can write $\rho(e_i)$ $(i=1,2,3,4)$ as follow:
\\$\left(\begin{array}{ccc}
a_{11}& a_{12} & a_{13}\\
   0  & a_{22} & a_{23}\\
   0  &   0    & a_{33}
\end{array}\right),$
$\left(\begin{array}{ccc}
b_{11}& b_{12} & b_{13}\\
   0  & b_{22} & b_{23}\\
   0  &   0    & b_{33}
\end{array}\right),$
$ \left(\begin{array}{ccc}
c_{11}& c_{12} & c_{13}\\
   0  & c_{22} & c_{23}\\
   0  &   0    & c_{33}
\end{array}\right),
$ $ \left(\begin{array}{ccc}
d_{11}& d_{12} & d_{13}\\
   0  & d_{22} & d_{23}\\
   0  &   0    & d_{33}
\end{array}\right).$

(1) $\frg=\frg_2(0,0)$. Then by the Lie brackets, we have
\begin{eqnarray*}
\hspace{1cm} (a)&& c_{11}=c_{22}=c_{33}=d_{11}=d_{22}=d_{33}=0;\\
(b)&&
c_{12}=b_{12}(a_{11}-a_{22})+a_{12}(b_{22}-b_{11}),c_{23}=b_{23}(a_{22}-a_{33})+a_{23}(b_{33}-b_{22}); \\
(c)&&c_{13}=b_{13}(a_{11}-a_{33})+a_{12}b_{23}+a_{13}b_{33}-b_{12}a_{23}-b_{11}a_{13};\\
(d)&&d_{12}=c_{12}(a_{11}-a_{22}),d_{23}=c_{23}(a_{22}-a_{33});\\
(e)&&d_{13}=c_{13}(a_{11}-a_{33})+a_{12}c_{23}-a_{23}c_{12};\\
(f)&&c_{12}(b_{11}-b_{22})=0, \quad d_{12}(b_{11}-b_{22})=0; \\
(g)&&c_{23}(b_{22}-b_{33})=0, \quad d_{23}(b_{22}-b_{33})=0;\\
(h)&&c_{13}(b_{11}-b_{33})+b_{12}c_{23}-b_{23}c_{12}=0;\\
(i)&&d_{13}(b_{11}-b_{33})+b_{12}d_{23}-b_{23}d_{12}=0;\\
(j)&&c_{12}(a_{11}-a_{22})^2=c_{12}(a_{11}-a_{22}),c_{23}(a_{22}-a_{33})^2=c_{23}(a_{22}-a_{33}).
\end{eqnarray*}
Case (1-i): $b_{11}$, $b_{22}$ and $b_{33}$ do not equal to each
other. Then $b_{11}\ne b_{22}$ or $b_{22}\ne b_{33}$. Without
losing the generality, we assume $b_{22}\ne b_{33}$ (in fact, the
discussion for the case $b_{11}\ne b_{22}$ is similar). Therefore
$c_{23}=d_{23}=0$ by $(g)$. Moreover, $c_{12}$ and $d_{12}$ cannot
be zero at the same time. Otherwise, $\rho (e_3)$ is proportional
to $\rho (e_4)$ by $(a)$, which is a contradiction. Hence we have
$b_{11}=b_{22}$ by $(f)$ and $c_{12}\ne 0$ by $(d)$. Thus, using
$(b)$, we know $c_{12}=b_{12}(a_{11}-a_{22})\neq 0$ which gives
$b_{12}\neq 0$ and $a_{11}\neq a_{22}$. By $(j)$, we have
$a_{11}-a_{22}=1$. On the other hand, by $(h)$ and $c_{23}=0$, we
know that $ c_{13}(b_{11}-b_{33})=c_{12}b_{23}$. If $c_{13}=0$,
then $b_{23}=0$ since $c_{12}\ne 0$. Therefore, by $(i)$ and
$d_{23}=0$, we have $d_{13}(b_{11}-b_{33})=0$. So $d_{13}=0$. It
gives that $\rho (e_3)$ is proportional to $\rho (e_4)$ by $(a)$,
which is a contradiction. If $c_{13}\neq 0$, then $b_{23}\neq 0$.
By $(d)$ and $a_{11}-a_{22}=1$, we have $c_{12}=d_{12}$. By $(i)$
and $d_{23}=0$, we have
$d_{13}(b_{11}-b_{33})=d_{12}b_{23}=c_{12}b_{23}=c_{13}(b_{11}-b_{33})$.
So $d_{13}=c_{13}$. Therefore $\rho(e_3)=\rho(e_4)$, which is a
contradiction.

\noindent Case (1-ii): $b_{11}=b_{22}=b_{33}$. Set
$e=\rho(e_2)-b_{11}I_3\ne 0 $ which is a strictly upper-triangle
matrix. Since $[e,\rho (e_3)]=[e,\rho (e_4)]=[\rho
(e_3),\rho(e_4)]=0,$, we have that $\{ e,\rho(e_3), \rho (e_4)\}$
spans an abelian subalgebra of the Lie algebra of the $3\times 3$
strictly upper-triangle matrices. Therefore $e$ is in the vector
space spanned by $\rho(e_3)$ and $\rho (e_4)$. We assume
$e=\alpha\rho(e_3)+\beta\rho(e_4)$. So
$\rho(e_3)=\rho([e_1,e_2])=[\rho(e_1),\rho(e_2)]=[\rho(e_1),e]=(\alpha+\beta)\rho(e_4)$,
which is  again a contradiction.

(2) $\frg=\frg_8(\frac{1}{4})$. With a similar discussion as in the
above case (1), we can know that $\mu (\frg_8(\frac{1}{4}))>3$.

On the other hand, for $\frg$ being $\frg_2(0,0)$ or
$\frg_8(\frac{1}{4})$, we can construct a 4-dimensional faithful
representation given in Table 3. \hfill $\Box$

%\begin{cor}
%With the notations given as above, if $\frg$ is isomorphic to one
%of Lie algebras appearing in Propositions 3.9 and 3.10, then
%$\mu(\frg)=4$.
%\end{cor}

%{\bf Proof}\quad \hfill $\Box$

At the end of this section, for every complex Lie algebra in
dimension $\leq 4$, we construct a faithful representations with the
lowest dimension and list it in the third column in Table 3. Let
$\{e_1,\cdots, e_m\}$ be a basis of $\frg$ and $e_{ij}$ be the
standard basis of the general linear Lie algebra $\frg\frl(n,\bbC)$.
We also give the corresponding left-symmetric algebras in the fourth
column in Table 3 which will be discussed in next section.

\vspace{1cm}

\mbox{}

\vspace{0.5cm}

{\small\begin{center} {\large T}{\rm ABLE 3}\end{center}
\begin{tabular}{|c|c|c|c|}
\hline $\frg$& $\mu(\frg)$ & Faithful representation &Left-symmetric algebra\\
\hline $\bbC$ & $1$ & $e_1\mapsto e_{11}$&$e_1*e_1=e_1$\\
\hline $\bbC^2$ & 2 & $e_1\mapsto e_{11}$,
$e_2\mapsto e_{22}$&$e_1*e_1=e_1$, $e_2*e_2=e_2$\\
\hline $\frr_2(\bbC)$ &2 & $e_1\mapsto e_{12}$,
$e_2\mapsto e_{22}$& $e_1*e_2=e_1$, $e_2*e_2=e_2$\\
\hline $\bbC^3$ &3 & $e_1\mapsto e_{11}$,
$e_2\mapsto e_{22}$, $e_3\mapsto e_{33}$& $e_1*e_1=e_1, e_2*e_2=e_2$, $e_3*e_3=e_3$\\
\hline $\frn_3(\bbC)$ &3 & $e_1\mapsto I_3+e_{12}-e_{23}$, &
$e_1*e_1=e_1-e_2+\frac{1}{2}e_3$,\\ & & $e_2\mapsto
e_{12}+e_{23}$, $e_3\mapsto 2e_{13}$ &
$e_1*e_2=e_2+\frac{1}{2}e_3$, $e_1*e_3=e_3$,
\\  &&&  $e_2*e_1=e_2-\frac{1}{2}e_3$, $e_2*e_2=\frac{1}{2}e_3$,
\\ &&&$e_3*e_1=e_3$\\
\hline $\frr_2(\bbC)\oplus\bbC$ & 2 & $e_1\mapsto
e_{12}$, $e_2\mapsto e_{22}$, $e_3\mapsto I_2$&\\
\hline $\frr_3(\bbC)$ & 3 & $e_1\mapsto e_{12}-e_{33}$, &
$e_1*e_1=-e_1+e_2$,  $e_1*e_3=e_2$,\\ && $e_2\mapsto e_{13}$,
$e_3\mapsto e_{23}$ &
 $e_2*e_1=-e_2$, $e_3*e_1=-e_3$\\
\hline $\frr_{3,\lambda}(\bbC)$ &3 & $e_1\mapsto
-e_{11}+(\lambda-1)e_{33}$, &$e_1*e_1=-e_1+(\lambda^2-\lambda)e_3$
\\ && $e_2\mapsto e_{21}$, $e_3\mapsto e_{31}$& $e_1*e_3=(\lambda-1)e_3$,
$e_2*e_1=-e_2$,\\ &&&$e_3*e_1=-e_3$\\
\hline $\frsl_2(\bbC)$ & 2 & $e_1\mapsto e_{12}$,
$e_2\mapsto e_{21}$, & \\&&$e_3\mapsto e_{11}-e_{22}$&\\
\hline $\bbC^4$ & 4&  $e_1\mapsto e_{11}$, $e_2\mapsto e_{22}$,
&$e_1*e_1=e_1$, $e_2*e_2=e_2$,
\\&&$e_3\mapsto e_{33}$, $e_4\mapsto
e_{44}$ & $e_3*e_3=e_3$, $e_4*e_4=e_4$\\
\hline $\frn_3(\bbC)\oplus\bbC$ & 3& $e_1\mapsto e_{12}$,
$e_2\mapsto e_{23}$,&\\&& $e_3\mapsto e_{13}$, $e_4\mapsto
I_3$& \\
\hline $\frr_2(\bbC)\oplus\bbC^2$ & 3 & $e_1\mapsto e_{13}$,
$e_2\mapsto e_{33}$,&\\&& $e_3\mapsto e_{22}$, $e_4\mapsto
I_3$&\\
\hline $\frr_3(\bbC)\oplus\bbC$ & 3 & $e_1\mapsto e_{12}-e_{33}$,
$e_2\mapsto e_{13}$,&\\&& $e_3\mapsto e_{23}$,
$e_4\mapsto I_3$&\\
\hline $\frr_{3,\lambda}(\bbC)\oplus\bbC$ & 3 & $e_1\mapsto
e_{11}+(1-\lambda)e_{33}$,&\\ && $e_2\mapsto e_{12}$,
$e_3\mapsto e_{13}$, $e_4\mapsto I_3$&\\
\hline $\frr_2(\bbC)\oplus\frr_2(\bbC)$ & 3 & $e_1\mapsto e_{13}$,
$e_2\mapsto e_{33}$, &\\&&$e_3\mapsto e_{12}$,
$e_4\mapsto e_{22}$&\\
\hline $\frsl_2(\bbC)\oplus\bbC$ & 2 & $e_1\mapsto e_{12}$,
$e_2\mapsto e_{21}$,&\\&& $e_3\mapsto e_{11}-e_{22}$, $e_4\mapsto
I_2$&\\\hline $\frn_4(\bbC)$ & 4 & $e_1\mapsto
I_4+e_{12}+e_{23}$,&$e_1*e_1=e_1+e_3+2e_4$, \\&& $e_2\mapsto
e_{34}$, $e_3\mapsto e_{24}$,
$e_4\mapsto e_{14}$& $e_1*e_2=e_2+e_3$, $e_1*e_3=e_3+e_4$, \\
&&&$e_1*e_4=e_4$, $e_2*e_1=e_2$, \\&&& $e_3*e_1=e_3$, $e_4*e_1=e_4$\\
\hline$\frg_1(\alpha)$ & 4 & $e_1\mapsto -e_{11}+(\alpha-1)e_{44}$,
&$e_1*e_1=-e_1+(\alpha^2-\alpha)e_4$,
\\&& $e_2\mapsto e_{21}$, $e_3\mapsto
e_{31}$, $e_4\mapsto e_{41}$& $e_1*e_4=(\alpha-1)e_4$, $e_2*e_1=-e_2$,\\
&&& $e_3*e_1=-e_3$, $e_4*e_1=-e_4$\\
\hline $\frg_2(\alpha,\beta) $ & 4 & $e_1\mapsto
-be_{11}-xe_{33}-ye_{44}+e_{32}$ &
$e_1*e_1=-be_1+(2b^3y-b^2y-2b^3+$\\&& $+e_{43}$,
$e_2\mapsto e_{21}$, $e_3\mapsto be_{21}+e_{31}$ & $b^2+bxy-bxy^2)e_2+(4b^2-3b+$\\
& & $e_4\mapsto
b^2e_{21}+(2b-x)e_{31}+e_{41}$,&$2b^2y-bxy-xy+xy^2-2by^2)e_3$
\\&& where $b\neq 0$ satisfies & $+(2-2b-by+y^2)e_4$,
\\ & & $b^3=\alpha-\beta
b+b^2$,&$e_1*e_2=e_3-be_2,e_1*e_3=e_4-be_3$,\\&&
$x^2-3bx+3b^2=-\beta+2b-x$,&$e_1*e_4=\alpha e_2-\beta
e_3+(1-b)e_4$,\\&& $y=3b-1-x$&$e_2*e_1=-be_2$,
 $e_3*e_1=-be_3$, \\&&&$e_4*e_1=-be_4$\\\hline
\end{tabular}}

\vspace{4cm}

\mbox{}

{\small\begin{center}{\large T}{\rm ABLE 3}
(Continued)\end{center}
\begin{tabular}{|c|c|c|c|}
\hline $\frg$& $\mu(\frg)$ & Faithful representation &Left-symmetric algebra\\
\hline  $\frg_3(\alpha)$ & 4 & $e_1\mapsto
-be_{11}-xe_{33}-ye_{44}+e_{32}$&
$e_1*e_1=-be_1+(2b^3y-2b^3+bxy-$\\&&
$+e_{43}$, $e_2\mapsto e_{21}$, $e_3\mapsto be_{21}+e_{31}$, & $bxy^2)e_2+(4b^2-b-4b^2y+bxy-$\\
& & $e_4\mapsto b^2e_{21}+(2b-x)e_{31}+e_{41}$,&
$xy+xy^2)e_3+(1-2b-by+y^2)e_4$,
\\&&   where $b^3=\alpha(b+1)$,&$e_1*e_2=e_3-be_2$,$e_1*e_3=e_4-be_3$
\\&& $x^2-3bx+3b^2-\alpha=0$, & $e_1*e_4=-be_4+\alpha(e_2+e_3)$\\&&$y=3b-x$ &$e_2*e_1=-be_2$, $e_3*e_1=-be_3$,\\&&&  $e_4*e_1=-be_4$\\
\hline $\frg_4$ & 4 & $e_1\mapsto
-e_{11}+e_{32}-\frac{3+\sqrt{3}i}{2}e_{33}$&
$e_1*e_1=-e_1+\frac{-1+\sqrt{3}i}{2}e_2+\frac{3+\sqrt{3}i}{2}e_3$
\\&&$+e_{43}-\frac{3-\sqrt{3}i}{2}e_{44}$, & $-(1+\sqrt{3}i)e_4$, $e_1*e_2=e_3-e_2$,
\\&&$e_2\mapsto e_{21}$, $e_3\mapsto e_{21}+e_{31}$, & $e_1*e_3=e_4-e_3$, $e_1*e_4=e_2-e_4$, \\
&& $e_4\mapsto e_{21}+\frac{1-\sqrt{3}i}{2}e_{31}+e_{41}$& $e_2*e_1=-e_2$, $e_3*e_1=-e_3$, \\
&&& $e_4*e_1=-e_4$\\
\hline $\frg_5$ & 4 & $e_1\mapsto
-\frac{1}{3}e_{44}+e_{23}$,&$e_1*e_1=-\frac{1}{3}e_1+\frac{1}{3}e_3$,\\&&
$e_2\mapsto e_{34}$, $e_3\mapsto e_{24},e_4\mapsto
e_{14}$&$e_1*e_2=e_3$, $e_2*e_1=-\frac{1}{3}e_2$,\\&&&
$e_3*e_1=-\frac{1}{3}e_3$, $e_4*e_1=-\frac{1}{3}e_4$\\
\hline $\frg_6$ & 3 & $e_1\mapsto e_{11}-e_{33}$, $e_2\mapsto
\frac{\sqrt{2}}{2}(e_{12}-e_{23})$,&\\&& $e_3\mapsto
\frac{\sqrt{2}}{2}(e_{12}+e_{23})$, $e_4\mapsto e_{13}$&\\
\hline $\frg_7$ &3& $e_1\mapsto e_{11}+e_{33}$, $e_2\mapsto
\frac{\sqrt{2}}{2}(e_{12}-e_{23})$,&\\&& $e_3\mapsto
\frac{\sqrt{2}}{2}(e_{12}+e_{23})$, $e_4\mapsto e_{13}$&\\
\hline $\frg_8(\alpha)$ & 3 & $e_1\mapsto
x_1e_{11}-x_2e_{33}$,&\\&& $e_2\mapsto
(1-4\alpha)^{-\frac{1}{4}}(e_{12}+e_{23})$, &\\
$(\alpha\neq\frac{1}{4})$ & & $e_3\mapsto
(1-4\alpha)^{-\frac{1}{4}}x_1e_{12}$&\\&& $
+(1-4\alpha)^{-\frac{1}{4}}x_2e_{23}$, $e_4\mapsto e_{14}$,
where&\\&& $x_1=\frac{1}{2}(1-(1-4\alpha)^{\frac{1}{2}})$,&\\ &&
$x_2=\frac{1}{2}(1+(1-4\alpha)^{\frac{1}{2}})$&\\
\hline $\frg_8(\frac{1}{4})$  & 4 & $e_1\mapsto
\frac{1}{2}(e_{11}-e_{44})+e_{23}$, &
$e_1*e_1=-\frac{1}{2}e_1-\frac{1}{4}e_2+\frac{1}{2}e_3+\frac{1}{2}e_4$,
\\&& $e_2\mapsto e_{12}+e_{13}+e_{24}+e_{34}$&$e_1*e_2=-\frac{1}{2}e_2+e_3+e_4$, \\ & &  $e_3\mapsto
\frac{{1}}{2}(e_{12}-e_{13}+3e_{24}+e_{34})$, & $e_1*e_3=-\frac{1}{4}e_2+\frac{1}{2}e_3+\frac{1}{4}e_4$,
\\&& $e_4\mapsto 2e_{14}$& $e_1*e_4=\frac{1}{2}e_4$, $e_2*e_1=-\frac{1}{2}e_2+e_4$ \\&&& $e_2*e_2=e_4$,
$e_2*e_3=e_4$,\\&&&$e_3*e_1=-\frac{1}{2}e_3+\frac{1}{4}e_4$\\&&& $e_3*e_3=\frac{1}{4}e_4$, $e_4*e_1=-\frac{1}{2}e_4$ \\
\hline
\end{tabular}}

\section{\'Etale affine representations and corresponding left-symmetric algebras}

Obviously, the construction of faithful representations at the end
of last section is not unique (in fact, there are a lot of faithful
representations). In this paper, we hope that the faithful
representation in the case $\dim\frg=\mu(\frg)$ can induce an
\'etale affine representation. In general, the construction is not
trivial (unfortunately at many cases we cannot use the construction
given in some proofs in section 3). For example, if the adjoint
representation of $\frg$ is faithful (hence $\frg$ is centerless),
it cannot induce an \'etale affine representation. Otherwise, in
this case, there exists an invertible derivation of $\frg$ (by
equation (2.5), $q^{-1}$ is a derivation) which implies $\frg$ is
nilpotent. Hence $\frg$ has a non-zero center, which is a
contradiction. Furthermore, we would also like to point out that
although the \'etale affine representations are not unique, either,
our construction has a natural and simple form  as follows (that is,
the ``construction rule").

\begin{thm}
For the faithful representations given in section 3 (as in Table 3),
if  $\mu (\frg)=\dim \frg$, then there exists an \'etale affine
representation $\Phi=(\rho, q):\frg\rightarrow \frgl(V)\oplus V$
with base zero, where $\rho:\frg\rightarrow \frgl (V)$ is just the
given faithful representation and
$$q(x)=\rho (x) (\sum_{i=1}^n e_i).\eqno (4.1)$$
\end{thm}

{\bf Proof}\quad Obviously, $q$ satisfies the equation (2.5). $q$ is
a linear isomorphism follows from the direct computation.\hfill
$\Box$

We also list the corresponding left-symmetric algebra structures
(with non-zero products) in the fourth column in Table 3 by
equation (2.6).

On the other hand, besides the semisimple Lie algebra
$\frs\frl(2,\bbC)$ (there does not exist a compatible left-symmetric
algebra structure [M]), for the Lie algebras with
$\dim\frg>\mu(\frg)$ in Table 3, we can also construct (in some
sense, like a kind of ``lifting") a faithful representation
$(\rho,V)$ with the same dimension of $\frg$ and it can induce an
\'etale affine representation $\Phi$ which still satisfies the
demands in Theorem 4.1. That is, $\Phi=(\rho, q):\frg\rightarrow
\frgl(V)\oplus V$, where $q(x)=\rho (x) (\sum_{i=1}^n e_i)$.

In Table 4, for those Lie algebras, we list $\dim\frg$ in column
2, faithful representation $\rho:\frg\rightarrow \frg\frl(\dim
\frg,\bbC)$ in column 3 and the corresponding left-symmetric
algebra (with non-zero products) in column 4.

\vspace{0.5cm}

{\small\begin{center} {\large T}{\rm ABLE 4}\end{center}
\begin{tabular}{|c|c|c|c|}
\hline $\frg$& $\dim\frg$ & Faithful representation &
Left-symmetric algebra
\\\hline
$\frr_2(\bbC)\oplus\bbC$ & 3& $e_1\mapsto e_{12}$, $e_2\mapsto
e_{22}$, $e_3\mapsto e_{33}$ & $e_1*e_2=e_1$, $e_2*e_2=e_2$,\\&&&
$e_3*e_3=e_3$\\
\hline $\frn_3(\bbC)\oplus\bbC$ & 4& $e_1\mapsto
e_{11}+e_{22}+e_{33}+e_{12}-e_{23}$,
&$e_1*e_1=e_1-e_2+\frac{1}{2}e_3$,\\&&$e_2\mapsto e_{12}+e_{23}$,
$e_3\mapsto 2e_{13}$,& $e_1*e_2=e_2+\frac{1}{2}e_3$,
$e_1*e_3=e_3$,\\&& $e_4\mapsto
e_{44}$&$e_2*e_1=e_2-\frac{1}{2}e_3$,
$e_2*e_2=\frac{1}{2}e_3$,\\&&& $e_3*e_1=e_3$,
 $e_4*e_4=e_4$\\\hline
$\frr_2(\bbC)\oplus\bbC^2$ & 4& $e_1\mapsto e_{12}$, $e_2\mapsto
e_{22}$, & $e_1*e_2=e_1$,
$e_2*e_2=e_2$,\\&& $e_3\mapsto e_{33}$, $e_4\mapsto e_{44}$&$e_3*e_3=e_3$, $e_4*e_4=e_4$\\
\hline $\frr_3(\bbC)\oplus\bbC$ & 4& $e_1\mapsto e_{12}-e_{33}$,
$e_2\mapsto e_{13}$, & $e_1*e_1=-e_1+e_2$,
$e_1*e_3=e_2$,\\&&$e_3\mapsto e_{23}$, $e_4\mapsto e_{44}$&
 $e_2*e_1=-e_2$, $e_3*e_1=-e_3$,\\ &&& $e_4*e_4=e_4$\\\hline
$\frr_{3,\lambda}(\bbC)\oplus\bbC$ & 4& $e_1\mapsto
-e_{11}+(\lambda-1)e_{33}$,&
$e_1*e_1=-e_1+(\lambda^2-\lambda)e_3$,\\&&$e_2\mapsto e_{21}$,
$e_3\mapsto e_{31}$, $e_4\mapsto e_{44}$&$e_1*e_3=(\lambda-1)e_3$,
 $e_2*e_1=-e_2$,\\&&& $e_3*e_1=-e_3$,  $e_4*e_4=e_4$\\\hline
$\frr_2(\bbC)\oplus\frr_2(\bbC)$ & 4& $e_1\mapsto e_{12}$,
$e_2\mapsto e_{22}$,&$e_1*e_2=e_1$, $e_2*e_2=e_2$,\\&& $e_3\mapsto
e_{34}$, $e_4\mapsto e_{44}$&$e_3*e_4=e_3$, $e_4*e_4=e_4$\\\hline
$\frsl_2(\bbC)\oplus\bbC$ & 4& $e_1\mapsto
2e_{12}+\frac{1}{2}e_{34}$, &$e_1*e_2=\frac{1}{2}(e_3+e_4)$,
$e_1*e_3=-e_1$, \\&&$e_2\mapsto
\frac{1}{2}e_{21}+2e_{43}$,&$e_1*e_4=e_1$,
$e_2*e_1=\frac{1}{2}(-e_3+e_4)$,\\&& $e_3\mapsto
e_{11}-e_{22}+e_{33}-e_{44}$,&$e_2*e_3=e_2$, $e_2*e_4=e_2$,\\&&
$e_4\mapsto I_4$& $e_3*e_1=e_1$, $e_3*e_2=-e_2$,\\
&&& $e_3*e_3=e_4$, $e_3*e_4=e_3$, \\&&&  $e_4*e_1=e_1$,
$e_4*e_2=e_2$, \\&&&    $e_4*e_3=e_3$, $e_4*e_4=e_4$\\\hline
$\frg_6$ & 4& $e_1\mapsto e_{11}-e_{44}$, &$e_1*e_1=2e_4-e_1$,
$e_1*e_2=2e_4$,\\&&$e_2\mapsto 2e_{12}+e_{34}$,& $e_1*e_3=e_4$,
$e_1*e_4=e_4$, \\&& $e_3\mapsto e_{13}+e_{24}$,&$e_2*e_1=2e_4-e_2$,
$e_2*e_3=2e_4$,\\&& $e_4\mapsto e_{14}$&$e_3*e_1=e_4-e_3$,
$e_3*e_2=e_4$,
\\&&&$e_4*e_1=-e_4$
\\\hline
\end{tabular}}

\vspace{2cm}

{\small\begin{center} {\large T}{\rm ABLE 4} (Continued)
\end{center}
\begin{tabular}{|c|c|c|c|}
\hline $\frg$& $\dim\frg$ & Faithful representation &
Left-symmetric algebra
\\
\hline $\frg_7$ & 4& $e_1\mapsto
e_{11}+2e_{33}+e_{44}$,&$e_1*e_1=e_1+e_2+e_3-e_4$, \\&&$e_2\mapsto
e_{12}-e_{24}+e_{34}$,&$e_1*e_2=e_2+e_3-\frac{1}{2}e_4$,\\&&
$e_3\mapsto
e_{12}+e_{24}+e_{34}$,&$e_1*e_3=e_2+e_3-\frac{1}{2}e_4$,\\&&
$e_4\mapsto 2e_{14}$&$e_1*e_4=e_4$,
$e_2*e_1=e_2-\frac{1}{2}e_4$,\\
&&&$e_2*e_2=-\frac{1}{2}e_4$, $e_2*e_3=\frac{1}{2}e_4$,\\&&&
$e_3*e_1=e_3-\frac{1}{2}e_4$,
$e_3*e_2=-\frac{1}{2}e_4$,\\&&&$e_3*e_3=\frac{1}{2}e_4$,
$e_4*e_1=e_4$\\\hline $\frg_8(\alpha)$ & 4 & $e_1\mapsto
(x-1)e_{22}-xe_{33}-e_{44}$,& $e_1*e_1=-e_1-\alpha
e_2+\frac{x}{2x-1}e_4$,\\ ($\alpha\neq\frac{1}{4}$)
 & & $e_2\mapsto
e_{12}+e_{24}+e_{34}$,&$e_1*e_2=-e_2+e_3+\frac{x}{2x-1}e_4$,\\&&
$e_3\mapsto (1-x)e_{12}+x_{24}$&$e_1*e_3=-\alpha
e_2+\frac{\alpha}{2x-1}e_4$,\\&& $+(1-x)e_{34}$, $e_4\mapsto
(2x-1)e_{14}$,&$e_2*e_1=-e_2+\frac{x}{2x-1}e_4$,\\&& where
$x=\frac{1}{2}(1+\sqrt{(1-4\alpha)})$ &
$e_2*e_2=\frac{1}{2x-1}e_4$,
$e_2*e_3=\frac{x}{2x-1}e_4$,\\&&&$e_3*e_1=-e_3+\frac{\alpha}{2x-1}e_4$,\\&&&$e_3*e_2=\frac{1-x}{2x-1}e_4$,
$e_3*e_3=\frac{\alpha}{2x-1}e_4$,\\&&&  $e_4*e_1=-e_4$
\\\hline
\end{tabular}}

%\begin{rmk}
%$\frg_8(\frac{1}{4})$
%\end{rmk}

\bigskip

\section*{Acknowledgements}

We thank Professor D. Burde for introducing us the history of
refining Ado's Theorem, bringing our attention to the reference
[GST] and the valuable discussion. We also thank the referee for the
important suggestion. This work was supported in part by the
National Natural Science Foundation of China (10501025, 10571091,
10621101), NKBRPC (2006CB805905), Program for New Century Excellent
Talents in University, Liu Hui Center for Applied Mathematics and
Youth Teachers Foundation of Tianjin University.

\end{document}